\documentclass{gtart}
\usepackage{amssymb,amsmath,amscd}

\input gtoutput
\volumenumber{2}\papernumber{7}\volumeyear{1998}
\pagenumbers{117}{144}\published{20 July 1998}
\proposed{David Gabai}\seconded{Walter Neumann, Wolfgang Metzler}
\received{21 December 1996}\revised{26 March 1998}
\accepted{17 July 1998}

\newtheorem{thm}{Theorem}[section]
\newtheorem{cor}[thm]{Corollary}

\newtheorem{lem}[thm]{Lemma}
\newtheorem{prop}[thm]{Proposition}
\newtheorem{dfn}[thm]{Definition}

\theoremstyle{remark}
\newtheorem{ex}[thm]{Example}
\newtheorem{rmk}[thm]{Remark}
\newtheorem{fact}[thm]{Fact}

\begin{document}

\title{Intersections in hyperbolic manifolds}
\author{Igor Belegradek}

\address{Department of Mathematics and Statistics\\McMaster University\\
1280 Main St West\\                                
Hamilton, ON  L8S 4K1\\                                    
Canada}

\email{belegi@icarus.math.mcmaster.ca}

\begin{abstract} We obtain some restrictions on the topology 
of infinite volume hyperbolic manifolds. In particular,
for any $n$ and any closed
negatively curved manifold $M$ of dimension $\ge 3$,
only finitely many hyperbolic $n$--manifolds
are total spaces of orientable vector bundles over $M$.
\end{abstract}

\asciiabstract{We obtain some restrictions on the topology 
of infinite volume hyperbolic manifolds. In particular,
for any n and any closed
negatively curved manifold M of dimension greater than 2,
only finitely many hyperbolic n-manifolds
are total spaces of orientable vector bundles over M.}

\primaryclass{30F40, 53C23, 57R20}
\secondaryclass{22E40, 32H20, 51M10}
\keywords{Hyperbolic manifold, intersection form, representation variety}

\maketitlepage

\section{Introduction}

A {\it hyperbolic manifold} is, by definition,  
a quotient of a negatively curved symmetric space by 
a discrete isometry group that acts freely. 
Recall that
negatively curved symmetric spaces
are hyperbolic spaces over the 
reals, complex numbers, quaternions,
or Cayley numbers.

The homotopy type of a hyperbolic manifold is determined
by its fundamental group.
Conversely,
for finite volume hyperbolic manifolds, fundamental group
determines the diffeomorphism type, thanks to the Mostow rigidity
theorem. 
By contrast infinite volume hyperbolic manifolds
with isomorphic fundamental groups can be very different 
topologically.
For example, the total spaces of many plane bundles
over closed surfaces carry hyperbolic metrics
(this and other examples are discussed in the 
section~\ref{S:Examples}).

In this paper we attempt to count hyperbolic manifolds 
up to intersection preserving homotopy equivalence.
A homotopy equivalence $f\co N\to L$ of oriented $n$--manifolds is called
{\it intersection preserving} if, for any $k$ and any pair of (singular) 
homology classes
$\alpha\in H_k(N)$ and $\beta\in H_{n-k}(N)$, their intersection number in
$N$ is equal to the intersection number of $f_*\alpha$ and $f_*\beta$
in $L$. For example, any map that is homotopic to an orientation-preserving
homeomorphism is an intersection preserving homotopy equivalence.
Conversely, oriented rank two vector bundles over a closed oriented surface
are isomorphic iff their total spaces are intersection preserving
homotopy equivalent.
\begin{thm}\label{intro:alg}
Let $\pi$ be the fundamental group of a finite 
aspherical cell complex and let $X$ be a negatively curved symmetric 
space. Let $\rho_k$ be a sequence of discrete injective
representations of $\pi$ into the group of orientation-preserving
isometries of $X$. 
Suppose that $\rho_k$ is precompact in the pointwise convergence topology.

Then the sequence of manifolds $X/\rho_k(\pi)$ falls into
finitely many intersection preserving homotopy equivalence classes.
\end{thm}

The space of conjugacy classes
of faithful discrete representations of a group 
$\pi$ into the isometry group of a negatively curved 
symmetric space is compact provided $\pi$
is finitely presented, not virtually
nilpotent and does not split over a virtually
nilpotent group (see~\ref{compactness thm}).
Since all the intersections in a hyperbolic manifold with
virtually nilpotent fundamental group are zero (see~\ref{nilpotent groups}),
we get the following.
\begin{cor}\label{no splitting implies
finiteness of intersection classes} 
Let $\pi$ be the fundamental group of a finite 
aspherical cell complex that does not split 
as an HNN--extension or a nontrivial amalgamated product over a
virtually nilpotent group. 

Then, for any $n$, the class of orientable
hyperbolic $n$--manifolds with fundamental
group isomorphic to $\pi$ breaks into
finitely many intersection preserving homotopy equivalence classes.
\end{cor}

A particular case of~\ref{no splitting implies
finiteness of intersection classes} was proved by Kapovich in~\cite{Kap2}.
Namely, he gave a proof for {\it real} hyperbolic $2m$--manifolds
homotopy equivalent to a closed orientable negatively
curved manifold of dimension $m\ge 3$.

If $\pi$ {\it does} split over a virtually nilpotent group,
the space of representations is usually noncompact.
For instance, this happens if $\pi$ is a surface group.
Yet, for real hyperbolic $4$--manifolds
homotopy equivalent to closed surfaces,
Kapovich~\cite{Kap2} proved a result
similar to~\ref{no splitting implies
finiteness of intersection classes}.

More generally, Kapovich~\cite{Kap2} proved that there is a universal
function $C(-,-)$ such that for any incompressible
singular surfaces $\Sigma_{g_1}$, $\Sigma_{g_2}$ in an 
oriented real hyperbolic $4$--manifold we have
$|\langle [\Sigma_{g_1}], [\Sigma_{g_2}]\rangle 
|\le C(g_1, g_2)$.
Reznikov \cite{Rez} showed that for any singular surfaces
$\Sigma_{g_1}$, $\Sigma_{g_2}$ in
a closed oriented negatively curved $4$--manifold $M$ with
the sectional curvature pinched between $-k^2$ and $-K^2$
there is a bound
$|\langle [\Sigma_{g_1}], [\Sigma_{g_2}]\rangle |\le 
C(g_1, g_2, k, K, \chi(M))$ for some universal function $C$.

Another way to classify hyperbolic manifolds is
up to tangential homotopy equivalence.
Recall that a homotopy equivalence of 
smooth manifolds $f\co N\to L$
is called tangential if the vector bundles
$f^*TL$ and $TN$ are stably isomorphic.
For example, any map that 
is homotopic to a diffeomorphism is a tangential 
homotopy equivalence. 
Conversely, a tangential homotopy equivalence
of open $n$--manifolds is homotopic to a diffeomorphism
provided $n>4$ and each of the manifolds
is the total space of a vector bundle
over a manifold of dimension $<n/2$~\cite[pp~226--228]{LS}.

An elementary argument 
(based on finiteness of the number of connected components
of representation varieties
and on finiteness of the number of symmetric spaces
of a given dimension) yields the following.
\begin{thm}\label{intr:tangent hom}
Let $\pi$ be a finitely presented group.
Then, for any $n$, the class of complete locally symmetric 
nonpositively curved
Riemannian $n$--manifolds with the fundamental group
isomorphic to $\pi$ falls into finitely many
tangential homotopy types.
\end{thm}

Knowing both the
intersection preserving and tangential homotopy
types sometimes
suffices to recover the manifold
up to finitely many possibilities.
\begin{thm}\label{intr:vector}
Let $M$ be a closed orientable negatively curved
manifold of dimension $\ge 3$
and $n>\dim(M)$ be an integer.
Let $f_k\co M\to N_k$ be a sequence of smooth
embeddings of $M$ into orientable hyperbolic $n$--manifolds
such that $f_k$ induces monomorphisms of fundamental groups.
Then the set of the normal bundles $\nu(f_k)$ of the
embeddings breaks into finitely many
isomorphism classes.

In particular,
only finitely many orientable hyperbolic $n$--manifolds
are total spaces of vector bundles over $M$.
\end{thm}
In some cases it is easy to decide when there exist
{\it infinitely} many rank $m$ vector bundles over
a given base. 
For example, by a simple K--theoretic argument
the set of isomorphism classes of
rank $m$ vector bundles over a finite cell complex 
$K$ is infinite provided $m\ge\dim (K)$
and $\oplus_k H^{4k}(K,\Bbb Q)\neq 0$.
Furthermore, 
oriented rank two vector bundles over $K$ are in one-to-one
correspondence with $H^2(K,\Bbb Z)$
via the Euler class. 
Note that many arithmetic closed real hyperbolic manifolds
have nonzero Betti numbers in all dimensions~\cite{MR}.
Any closed complex hyperbolic manifold has nonzero 
even Betti numbers because the powers
of the K\"ahler form are noncohomologous to zero.
Similarly, for each $k$, closed quaternion hyperbolic
manifolds have nonzero $4k$th Betti numbers. 

Examples of vector bundles with hyperbolic total spaces
are given in the section~\ref{S:Examples}.
Note that according to a result of Anderson~\cite{And}
the total space of any vector bundle over a 
closed negatively curved manifold
admits a complete metric with the sectional
curvature pinched between two negative constants.

I am grateful to 
Bill Goldman,  Misha Kapovich, and Jonathan Rosenberg
for many helpful discussions.
I would like to thank Heinz Helling and SFB-343 at the University of Bielefeld
for support and hospitality.
\paragraph{Outline of the paper.}
The section~\ref{S:Examples} is a collection of 
examples of hyperbolic manifolds.
Some invariants of maps and representations are defined in
sections~\ref{S:Invariants maps} and \ref{S:Invariants of representations}.
Sections~\ref{S:spaces repr}, \ref{S:stable range}, 
and~\ref{S:loc symm} are devoted to a proof of
the theorem~\ref{intr:tangent hom} and other related results.
Sections~\ref{S:geom conv}, \ref{S:alg conv}, and~\ref{S:comp thm} 
contain background on algebraic 
and geometric convergence needed
for the main theorem which is proved in section~\ref{S:main thm}.
Theorems~\ref{intro:alg} and~\ref{no splitting implies
finiteness of intersection classes} are proved 
in section~\ref{S:intersection}. 
Finally, theorem~\ref{intr:vector}
is proved in section~\ref{S:vector}.

\section{Examples}
\label{S:Examples}
To help the reader appreciate the results stated above
we collect some relevant examples
of hyperbolic manifolds. 

\begin{ex}\bf(Plane bundles over closed surfaces)\rm\stdspace
Total spaces of rank two vector bundles over closed
surfaces often admit hyperbolic metrics.
For instance, an orientable 
$\Bbb R^2$--bundle over a closed oriented 
surface of genus $g$ admits a real hyperbolic structure provided 
the Euler number $e$ of the bundle satisfies $|e|<g$~\cite{Luo1, Luo2} 
(cf~\cite{GLT},~\cite{Kap1},~\cite{Kui1, Kui2}).
Complex hyperbolic structures exist on orientable 
$\Bbb R^2$--bundles over 
closed oriented surfaces when $|e+2g-2|<g$~\cite{GKL}.

Note that the Euler number is equal to
the self-intersection number
of the zero section, hence the total spaces
of bundles with different Euler classes
are not intersection homotopy equivalent.

For nonorientable bundles over nonorientable surfaces 
the condition $|e|\le [\frac{g}{8}]$ on the twisted
Euler number implies the existence
of a real hyperbolic structure~\cite{Bel1}.
\end{ex}

\begin{ex}\bf(Plane bundles over closed 
hyperbolic $3$--manifolds)\rm\nl
Total spaces of plane bundles over closed 
hyperbolic $3$--manifolds sometimes carry
hyperbolic metrics~\cite{Bel2}.

In fact, it can be deduced from~\cite{Bel2} that
for every $k$ there exists a closed oriented
real hyperbolic $3$--manifold $M=M(k)$ and 
oriented real hyperbolic $5$--manifolds
$N_1,\dots , N_k$ that are total spaces of
plane bundles over $M$ and such that no two of them 
are intersection preserving
homotopy equivalent. 
\end{ex}

\begin{ex}\bf(Fundamental group at infinity)\rm\stdspace
There are real hyperbolic $4$--manifolds
that are intersection homotopy equivalent but not homeomorphic
to plane bundles over closed surfaces~\cite{GLT}.  
The invariant that distinguishes these manifolds
from vector bundles is the fundamental group at infinity.

Even more surprising examples
were given in~\cite{Mat} and~\cite{GK}.
Namely, there are orientable 
real hyperbolic $4$--manifolds that are
homotopy equivalent but not homeomorphic to handlebodies;
these manifolds have nontrivial fundamental group at
infinity. 
Note that if $N$ and $L$ are orientable $4$--manifolds
that are homotopy equivalent to a handlebody,
then each homotopy equivalence $N\to L$
is both tangential and intersection preserving.
\end{ex}

\begin{ex}\bf(Nonorientable line bundles)\rm\stdspace
Here is a simple way to produce homotopy equivalent
hyperbolic manifolds that are not 
tangentially homotopy equivalent.

First, note that via the inclusion
$\bold{O}(n,1)\hookrightarrow\bold{O}(n+1,1)$, 
any discrete subgroup of $\bold{O}(n,1)$
can be thought of as a discrete subgroup of $\bold{O}(n+1,1)$
that stabilizes a subspace 
$\bold{H}^n_\Bbb R\subset\bold{H}^{n+1}_\Bbb R$.
The orthogonal projection $\bold{H}^{n+1}_\Bbb R\to\bold{H}^{n}_\Bbb R$
is an $\bold{O}(n,1)$--equivariant line bundle over
$\bold{H}^{n}_\Bbb R$.
In particular, given a real hyperbolic manifold $M$,
the manifold $M\times\Bbb R$ carries a real hyperbolic metric.

If $H^1(M,\Bbb Z_2)\neq 0$,
this construction can be twisted to produce
nonorientable line bundles over $M=\bold{H}^n_\Bbb R/\pi_1(M)$ with real 
hyperbolic metrics.
Indeed, a nonzero element 
$w\in H^1(M,\Bbb Z_2)\cong\mathrm{Hom}(\pi_1(M),\Bbb Z_2)$
defines an epimorphism $w\co \pi_1(M)\to\Bbb Z_2$.
Make $\Bbb Z_2$ act on $\bold{H}^{n+1}_\Bbb R$
as the reflection in $\bold{H}^{n}_\Bbb R$.
Then let $\pi_1(M)$ act on line bundle
$\bold{H}^{n}_\Bbb R\times\Bbb R\cong\bold{H}^{n+1}_\Bbb R$ 
by $\gamma(x,t)=(\gamma(x),w(\gamma)(t))$.
The quotient $\bold{H}^{n+1}_\Bbb R/\pi_1(M)$ is the total space
of a line bundle over $M$ with the first Stiefel--Whitney class $w$.

In particular, line bundles with different
first Stiefel--Whitney classes have total space that are
homotopy equivalent but {\it not} tangentially homotopy equivalent.  
(Tangential homotopy equivalences 
preserve Stiefel--Whitney classes yet
$w_1(\bold{H}^{n+1}_\Bbb R/\pi_1(M))=w_1(M)+w$.)
Thus, we get many tangentially homotopy inequivalent
manifolds in a given homotopy type.   
\end{ex}

\section{Invariants of continuous maps}
\label{S:Invariants maps}
\begin{dfn} 
Let $B$ be a topological space and $S_B$ be a set.
Let $\iota$ be a map that, given a smooth manifold $N$, 
and a continuous map from $B$ into $N$,
produces an element of $S_B$. 
We call $\iota$ an {\it invariant of maps of $B$}
if the two following conditions hold:

{\rm(1)}\stdspace Homotopic maps $f_1\co B\to N$ and 
$f_2\co B\to N$ have the same invariant.

{\rm(2)}\stdspace Let $h\co  N\to L$ be a
diffeomorphism of $N$ onto an open subset of $L$.
Then, for any continuous map $f\co B\to N$, 
the maps $f\co B\to N$ and $h\circ f\co B\to L$
have the same invariant.
\end{dfn}

There is a version of this definition for maps into
oriented manifolds.
Namely, we require that the target manifold is oriented and 
the diffeomorphism $h$ preserves orientation.
In that case we say that $\iota$ is an 
{\it invariant of maps of\/ $B$ into oriented manifolds}. 

\begin{ex}\bf(Tangent bundle)\rm\stdspace
\label{Tangent bundle}
Assume $B$ is paracompact and $S_B$ is the set 
of isomorphism classes of real
vector bundles over $B$.
Given a continuous map 
$f\co B\to N$, set $\tau (f\co B\to N)=f^\#TN$,
the isomorphism class of
the pullback of the tangent bundle to $N$
under $f$.
Clearly, $\tau$ is an invariant.
\end{ex}

\begin{ex}\bf(Intersection number in oriented $n$--manifolds)\rm\stdspace 
\label{Intersection number in oriented $n$--manifolds}
Assume $B$ is compact.
Fix two cohomology classes 
$\alpha\in H_{m}(B)$ and $\beta\in H_{n-m}(B)$.
(In this paper we always use singular (co)homology
with integer coefficients unless stated otherwise.)

Let $f\co B\to N$ be a continuous map of a
compact topological space $B$ into
an oriented $n$--manifold $N$ where $\mathrm{dim}(N)=n$.
Set $I_{n,\alpha,\beta}(f)$ to be the 
intersection number $I(f_*\alpha, f_*\beta)$
of $f_*\alpha$ and $f_*\beta$ in $N$. 
We next show that $I_{n,\alpha,\beta}$
is an integer-valued invariant of maps into 
oriented manifolds.

Recall that the intersection number $I(f_*\alpha, f_*\beta)$
can be defined as follows.
Start with an arbitrary compact subset $K$ of $N$ that contains $f(B)$.
Let $A\in H^{n-m}(N,N\backslash K)$ and 
$B\in H^{m}(N,N\backslash K)$ be the Poincar\'e
duals of $f_*\alpha\in H_{m}(K)$ and $f_*\beta\in H_{n-m}(K)$, 
respectively.
Then, set $I(f_*\alpha, f_*\beta)=\langle A\cup B, [N, N\backslash K]\rangle$
where $[N, N\backslash K]$ is the fundamental class of $N$ near 
$K$~\cite[VII.13.5]{Dol}. 
Note that
$$I(f_*\alpha, f_*\beta)=\langle A\cup B, [N, N\backslash K]\rangle=
\langle A, B\cap [N, N\backslash K]\rangle=\langle A, f_*\beta\rangle.$$
The following commutative diagram shows that $I(f_*\alpha, f_*\beta)$ 
is independent of $K$.
$$
\CD 
H_m(B) @>f_*>> H_m(N,N\backslash f(B)) @>D^{-1}>> 
H^{n-m}(f(B)) @>f^*>> H^{n-m}(B)\\
@Vf^*VV @Ai_*AA @AAi^*A @| \\
H_m(K) @>i_*>> H_m(N,N\backslash K)@>D^{-1}>> 
H^{n-m}(K) @>f^*>> H^{n-m}(B)\\
\endCD
$$

Note that $I_{n,\alpha,\beta}$ is an invariant.
Indeed, property $\mathrm{(1)}$ holds trivially;
property $\mathrm{(2)}$ is verified in~\cite[VIII.13.21(c)]{Dol}.
\end{ex}

We say that an invariant of maps is {\it liftable}
if in part $\mathrm{(2)}$ of the definition
the word ``diffeomorphism'' can be replaced by a
``covering map''.
For example, tangent bundle is a liftable invariant.
Intersection number is not liftable.
The following proposition shows to what extent 
it can be repaired.
 
\begin{prop}\label{liftable prop} Let $p\co \tilde{N}\to N$ be a
covering map of manifolds and let $B$ 
be a finite connected CW--complex.
Suppose that $f\co B\to\tilde{N}$ is a 
map such that $p\circ f\co B\to N$ is an embedding 
(ie a homeomorphism onto its image).

Then $\iota(f)=\iota(p\circ f)$
for any invariant of maps $\iota$.
\end{prop}

\begin{proof} 
Since $p\circ f\co B\to L$ is an embedding, so is $f$. 
Then, the map $p|_{f(B)}\co$ $f(B)\to p(f(B))$ is a homeomorphism.
Using compactness of $f(B)$, one can
find an open neighborhood $U$ of $f(B)$ such that
$p|_{U}\co U\to p(U)$ is a diffeomorphism.
Since invariants $\iota(f)$ and 
$\iota (p\circ f)$ can be computed in $U$ 
and $p(U)$, respectively, we conclude 
$\iota(f)=\iota(p\circ f)$.
\end{proof}

\section{Invariants of representations}
\label{S:Invariants of representations}
Assume $X$ is a smooth contractible manifold and 
let $\mathrm{Diffeo}(X)$ be the 
group of all self-diffeomorphisms of $X$ equipped with
compact-open topology.
Let $\pi$ be the fundamental group of a 
finite-dimensional CW--complex $K$
with universal cover $\tilde{K}$.

We refer to a group homomorphism
$\rho\co \pi\to\mathrm{Diffeo}(X)$ as a
{\it representation}.
To any representation $\rho\co \pi\to\mathrm{Diffeo}(X)$,
we associate a 
continuous $\rho$--equivariant map $\tilde{K}\to X$ as follows.
Consider the $X$--bundle
$\tilde{K}\times_\rho X$ over $K$
where $\tilde{K}\times_\rho X$
is the quotient of $\tilde{K}\times X$
by the following action of $\pi$ 
$$\gamma(\tilde k,x)=(\gamma(\tilde k), \rho(\gamma)(x)),
\ \ \gamma\in\pi.$$ 
Since $X$ is contractible, the bundle has a section
that is unique up to homotopy through sections.
Any section can be lifted to
a $\rho$--equivariant continuous map 
$\tilde{K}\to \tilde{K}\times X$.
Projecting to $X$, we get a $\rho$--equivariant continuous map 
$\tilde{K}\to X$.

Note that any two $\rho$--equivariant continuous maps 
$\tilde{g}, \tilde{f}\co \tilde{K}\to X$, 
are $\rho$--equivariantly homotopic.
(Indeed, $\tilde{f}$ and $\tilde{g}$ descend to sections
$K\to \tilde{K}\times_\rho X$ that must be homotopic.
This homotopy lifts to a $\rho$--equivariant homotopy
of $\tilde{f}$ and $\tilde{g}$.) 
 
Assume now that $\rho(\pi)$ acts freely and 
properly discontinuously on $X$.
Then the map $\tilde{f}$ descends to a continuous map 
$f\co K\to X/\rho(\pi_1(K))$.
We say that $\rho$ is {\it induced} by $f$.

Let $\iota$ be an invariant of continuous maps of $K$.
Given a representation $\rho$ such that
$\rho(\pi)$ acts freely and 
properly discontinuously on $X$, set
$\iota(\rho)$ to be $\iota(f)$
where $\rho$ is induced by $f$.
We say $\rho$ is an {\it invariant of representations of $\pi_1(K)$}. 

Similarly, any invariant $\iota$ of continuous maps of $K$
into oriented manifolds 
defines an invariant of representations
into the group of orientation-preserving
diffeomorphisms of $X$.  

Note that representations conjugate 
by a diffeomorphism $\phi$ of $X$ have the
same invariants.
(Indeed, if $\tilde f\co\tilde{K}\to X$ is a 
$\rho$--equivariant map, the map $\phi\circ\tilde f$
is $\phi\circ\rho\circ\phi^{-1}$--equivariant.)
The same is true for invariants of orientation-preserving
representations when $\phi$ is orientation-preserving.

\begin{ex}\bf(Tangent bundle)
\label{invariant tau}\rm\stdspace
Let $\tau$ be the invariant of maps defined in~\ref{Tangent bundle}.
Then, for any representation $\rho$ such that
$\rho(\pi)$ acts freely and 
properly discontinuously on $X$, 
let $\tau(\rho)$ be the pullback of the tangent 
bundle to $X/\rho(\pi)$ via a map $f\co K\to X/\rho(\pi)$ 
that induces $\rho$.

In fact, $\tau$ can be defined for any representation
as follows.
Look at the ``vertical'' bundle $\tilde{K}\times_\rho TX$ over 
$\tilde{K}\times_\rho X$ where $TX$ is the tangent bundle to $X$.
Set $\tau(\rho)$ to be the pullback of the vertical bundle
via a section $K\to\tilde{K}\times_\rho X$. 
Thus, to every representation $\rho\co \pi_1(K)\to\mathrm{Diffeo}(X)$ 
we associated a vector bundle $\tau(\rho)$ of rank $\dim(X)$ over $K$.
\end{ex}

\begin{ex}\bf(Intersection \kern-0.7pt number for
orientation \kern-0.7pt  preserving actions)\rm\stdspace
Assume the cell complex $K$ is finite
and choose an orientation on $X$ (which
makes sense because, like any contractible
manifold, $X$ is orientable).

Fix two cohomology classes 
$\alpha\in H_{m}(K)$ and $\beta\in H_{n-m}(K)$.
Let $I_{n,\alpha,\beta}$ is an invariant
of maps defined in~\ref{Intersection number in oriented $n$--manifolds}
where $n=\mathrm{dim}(X)$.

Let $\rho$ be a representation of $\pi$
into the group of orientation-preserving
diffeomorphisms of $X$ such that
$\rho(\pi)$ acts freely and 
properly discontinuously on $X$.
Then let $I_{n,\alpha,\beta}(\rho)$ be the 
intersection number $I(f_*\alpha, f_*\beta)$
of $f_*\alpha$ and $f_*\beta$ in $X/\rho(\pi)$ 
where $f\co K\to X/\rho(\pi)$ is a map
that induces $\rho$.
\end{ex}

\section{Spaces of representations and
tangential homotopy equivalence}
\label{S:spaces repr}
Let $X$ be a smooth contractible manifold and 
let $\mathrm{Diffeo}(X)$ be the 
group of all self-diffeomorphisms of $X$ equipped with
compact--open topology.
We equip the space of representations 
$\mathrm{Hom}(\pi,\mathrm{Diffeo}(X))$ with 
the pointwise convergence topology, ie
a sequence of representations $\rho_k$ converges
to $\rho$ provided $\rho_k(\gamma)$ converges to $\rho(\gamma)$ 
for each $\gamma\in\pi$.

In the next two sections we explore the
consequences of the following observation.

\begin{prop}\label{first prop in
Spaces of representations and
tangential homotopy equivalence} 
Let $\rho_0$ and $\rho_1$ be injective representations
of $\pi$ into\break $\mathrm{Diffeo}(X)$ such that
the groups $\rho_0(\pi)$ and $\rho_1(\pi)$
act freely and properly discontinuously on $X$. 
Suppose that $\rho_0$ and $\rho_1$
can be joined by a continuous path of representations
$\rho_t\co \pi\to\mathrm{Diffeo}(X)$ 
(where continuous means that, for every $\gamma\in\pi$,
the map $\rho_t(\gamma)\co [0,1]\to\mathrm{Diffeo}(X)$
is continuous).

Then the homotopy equivalence of manifolds
$X/\rho_0(\pi)$ and $X/\rho_1(\pi)$ induced by
$\rho_1\circ(\rho_0)^{-1}$ is tangential. 
\end{prop}
\begin{proof}
Since $\rho_t\co \pi\to\mathrm{Diffeo}(X)$
is a continuous path of representations,
the covering homotopy theorem implies that
the bundles $\tau(\rho_0)$ and $\tau(\rho_1)$
are isomorphic (the invariant $\tau$
is defined in~\ref{invariant tau}). 

Let $f_0\co K\to X/\rho_0(\pi)$ and
$f_1\co K\to X/\rho_1(\pi)$ be homotopy equivalences
that induce $\rho_0$ and $\rho_1$,
respectively. 
For $i=1,2$ the bundle $\tau(\rho_i)$ is
isomorphic to the pullback of the tangent bundle to
$X/\rho_i(\pi)$ via $f_i$.
Thus $f_1\circ (f_0)^{-1}$ is 
a tangential homotopy equivalence.
\end{proof}

\begin{rmk}
In fact, the covering homotopy theorem implies that
$\tau$ is constant on any path-connected component
of the space $\mathrm{Hom}(\pi_1(K),\mathrm{Diffeo}(X))$.
\end{rmk}

\begin{cor} \label{smooth conjugacy}
Under the assumptions 
of\/ {\rm\ref{first prop in
Spaces of representations and
tangential homotopy equivalence}}, suppose that
$\dim(X)\ge 5$ and that
each of the manifolds $X/\rho_0(\pi)$ 
and $X/\rho_1(\pi)$ is homeomorphic to the total space 
of a vector bundle over a manifold of dimension $<\dim(X)/2$.

Then $\rho_0$ and $\rho_1$ are smoothly conjugate on $X$.
\end{cor}
\begin{proof}
According to~\cite[pp 226--228]{LS}, the homotopy equivalence
$f_1\circ (f_0)^{-1}$  is homotopic to a diffeomorphism.
Hence, $\rho_0$ and $\rho_1$ are smoothly conjugate on $X$.
\end{proof}

\begin{rmk}
More precisely, the result proved in~\cite[pp 226--228]{LS} is as follows.
Let $f\co N_0\to N_1$ be a tangential homotopy equivalence of 
smooth $n$--manifolds with $n\ge 5$.
If each of the manifolds is 
homeomorphic to the interior
of a regular neighbourhood of a simplicial complex
of dimension $<n/2$,
then $f$ is homotopic to
a diffeomorphism. 
\end{rmk}

\section{Discrete representations in stable range}
\label{S:stable range}
Suppose $X$ is a symmetric space of nonpositive sectional curvature. 
Note that for any discrete torsion-free subgroup $\Gamma\le\mathrm{Isom}(X)$
that stabilizes a totally geodesic submanifold $Y$,
the exponential map identifies the quotient $X/\Gamma$ and 
the total space of normal bundle of $Y/\Gamma$ in $X/\Gamma$.
Applying~\ref{first prop in
Spaces of representations and
tangential homotopy equivalence}, we deduce the following.

\begin{thm}
Let $X$ be a nonpositively curved symmetric space of dimension
$\ge 5$ and let $\pi$ be a group.
Let $\rho_1$ and $\rho_2$ be injective discrete representations
of $\pi$ into the isometry group of $X$
that lie in the same path-connected component of the space
$\mathrm{Hom}(\pi,\mathrm{Isom}(X))$.
Suppose that each of the representation $\rho_1$ and $\rho_2$
stabilizes a totally geodesic subspace of dimension
$<\dim(X)/2$.

Then $\rho_1$ and $\rho_2$ are smoothly conjugate on $X$.
\end{thm}

\begin{rmk} 
Of course, the above argument works in 
other geometries as well.
Here is a sample result for complete
affine manifolds.

Let $\rho_0$ and $\rho_1$ be injective representations
of a group $\pi$ into $\mathrm{Aff}(\Bbb R^n)$ 
that lie in the same path-connected component of the 
space of representations $\mathrm{Hom}(\pi,$ $\mathrm{Aff}(\Bbb R^n))$.
Assume that the groups $\rho_0(\pi)$ and $\rho_1(\pi)$
act freely and properly discontinuously on $\Bbb R^n$
and, furthermore, suppose that $\rho_0(\pi)$ and $\rho_1(\pi)$
are contained in $\mathrm{Aff}(\Bbb R^k)\subset\mathrm{Aff}(\Bbb R^n)$.

Since the coordinate projection $\Bbb R^n\to\Bbb R^k$
is $\mathrm{Aff}(\Bbb R^k)$--equivariant,
the manifolds $\Bbb R^n/\rho_0(\pi)$ and $\Bbb R^n/\rho_1(\pi)$
are the total spaces of vector bundles over the manifolds
$\Bbb R^k/\rho_0(\pi)$ and $\Bbb R^k/\rho_1(\pi)$,
respectively.
Then~\ref{smooth conjugacy} implies that
$\rho_0$ and $\rho_1$ 
are conjugate by a diffeomorphism of $\Bbb R^n$
provided $n\ge 5$ and $k<n/2$.
\end{rmk}

\section{Locally symmetric nonpositively curved manifolds up to 
tangential homotopy equivalence}
\label{S:loc symm}
Let $G$ be a subgroup of $\mathrm{Diffeo}(X)$
such that the space of representations $\mathrm{Hom}(\pi, G)$
has finitely many path-connected components.
Then~\ref{first prop in
Spaces of representations and
tangential homotopy equivalence} implies that 
the class of manifolds of the form $X/\rho(\pi)$ where
$\rho$ is injective and $\rho(\pi)$ acts
freely and properly discontinuously on $X$
falls into finitely many tangential equivalence classes.

For example, $\mathrm{Hom}(\pi, G)$
has finitely many path-connected components if
$\pi$ is finitely presented and
$G$ is either real algebraic, or complex algebraic,
or semisimple with finite center~\cite{Gol}.
In particular, the following is true.

\begin{thm}\label{thm on finitely many tan hom types}
Let $\pi$ be a finitely presented torsion-free group
and let $X$ be a nonpositively curved symmetric space.
Then the class of manifolds of the form $X/\rho(\pi)$,
where $\rho\in\mathrm{Hom}(\pi, \mathrm{Isom}(X))$ is a faithful 
discrete representation, falls into finitely many
tangential homotopy types.
\end{thm}
\begin{proof} 
Represent $X$ as a Riemannian product $Y\times\Bbb R^k$
where $Y$ is a nonpositively curved symmetric space
without Euclidean factors. 
By de Rham's theorem this decomposition is unique,
so $\mathrm{Isom}(X)\cong\mathrm{Isom}(Y)\times\mathrm{Isom}(\Bbb R^k)$.
The group $\mathrm{Isom}(Y)$ is semisimple with trivial center,
hence the analytic variety $\mathrm{Hom}(\pi,\mathrm{Isom}(Y))$
has finitely many path-connected components~\cite[p 567]{Gol}.
The same is true for $\mathrm{Hom}(\pi,\mathrm{Isom}(\Bbb R^k))$
because $\mathrm{Isom}(\Bbb R^k)$ is real algebraic~\cite[p 567]{Gol}.
Hence the analytic variety
$$\mathrm{Hom}(\pi,\mathrm{Isom}(X))\cong
\mathrm{Hom}(\pi,\mathrm{Isom}(Y))\times
\mathrm{Hom}(\pi,\mathrm{Isom}(\Bbb R^k))$$
has finitely many connected components.
\end{proof}

\begin{cor}
Let $\pi$ be a finitely presented torsion-free group.
Then, for any $n$,
the class of locally symmetric complete nonpositively curved
$n$--manifolds with the fundamental group isomorphic to $\pi$
falls into into finitely many
tangential homotopy types.
\end{cor}
\begin{proof}
For any $n$ there exist only finitely many 
nonpositively curved symmetric spaces of dimension $n$.
Hence~\ref{thm on finitely many tan hom types} applies.
\end{proof}

\section{Geometric convergence}
\label{S:geom conv}
In this section we discuss some basic facts on geometric
convergence.
The notion of geometric convergence was introduced
by Chaubaty (see~\cite{Bou}).
More details relevant to our exposition
can be found in~\cite{CEG},~\cite{Lok} and~\cite{BP}.
In this section we let $G$ be a Lie group
equipped with some left invariant Riemmanian metric.

\begin{dfn} Let
$C(G)$ be the set of all closed subgroups of $G$.
Define a topology on $C(G)$ as follows.
We say that a sequence $\{\Gamma_n\}\in C(G)$ converges
to $\Gamma_{\mathrm{geo}}\in C(G)$ geometrically if
the following two conditions hold:

{\rm(1)}\stdspace If $\gamma_{n_{k}}\in \Gamma_{n_{k}}$ converges to
$\gamma \in G$, then $\gamma \in \Gamma_{\mathrm{geo}}$.

{\rm(2)}\stdspace If $\gamma \in \Gamma_{\mathrm{geo}}$, then there
is a sequence $\gamma_n\in \Gamma_n$ with
$\gamma_n\to \gamma$ in $G$.
\end{dfn}

\begin{fact}(\cite{CEG})\stdspace
The space $C(G)$ is compact and metrizable.
\end{fact}

\begin{fact}(\cite{Lok})\stdspace
$\Gamma_n\to \Gamma_{\mathrm{geo}}$ iff
 for every compact subset $K\subset G$ the sequence
$\Gamma_n\cap K\to \Gamma_{\mathrm{geo}}\cap K$ 
in the Hausdorff topology
(ie for any $\varepsilon >0$, there is $N$ such that, 
if $n>N$,
then $\Gamma_n\cap K$ lies in the $\varepsilon$--neihgborhood of
$\Gamma_{\mathrm{geo}}\cap K$ and $\Gamma_{\mathrm{geo}}\cap K$
lies in the $\varepsilon$--neihgborhood of $\Gamma_n\cap K$).
\end{fact}

\begin{fact}\label{lok: discrete groups are open}
(\cite{Lok})\stdspace
Let $\Gamma_{\mathrm{geo}}\subset G$ is a discrete
subgroup. Then there is $\varepsilon >0$ such that, for any
$\Gamma_n\to \Gamma_{\mathrm{geo}}$ in $C(G)$ and
any compact $K\subset G$, there is $N$ such that,
if $n>N$ and $\gamma\in \Gamma_{\mathrm{geo}}\cap K$, then
there is a unique $\gamma_n\in \Gamma_n$
that is $\varepsilon$--close to
$\gamma$.

In particular, $\Gamma_n$ is discrete for $n>N$, since $e\in \Gamma_n$
is the only element of $\Gamma_n$,
that is in the $\varepsilon$--neighbourhood
of the identity.\end{fact}

\begin{rmk}\label{limit is torsion free}
Let $\Gamma_n $ be a sequence of torsion-free groups
converging to a discrete group $\Gamma_{\mathrm{geo}}$ in $C(G)$.
Then $\Gamma_{\mathrm{geo}}$ is torsion-free.
Indeed, choose $\gamma\in \Gamma_{\mathrm{geo}}$ with $\gamma^k=e$.
Find a sequence $\gamma_n\to \gamma$. Then $\gamma_n^k\to e$. 
By~\ref{lok: discrete groups are open}
we have $\gamma_n^k=e$ for large $n$. Since $\Gamma_n$ are torsion-free,
$k=1$ as desired.
\end{rmk}

\begin{thm}[\cite{Lok}] 
\label{lok: lipschitz convergence}
Let $X$ be a simply connected homogeneous
Riemannian manifold and $G$ be a transitive group of isometries.
Let $\Gamma_n$ be a sequence of torsion-free subgroups of $G$
converging geometrically to a discrete subgroup $\Gamma_{\mathrm{geo}}$.
Let $U\subset X$ be a relatively compact open set.

Then, if $n$ is
large enough, there exists a relatively
compact open set $V$ with $U\subset V$
and smooth embeddings $\tilde\varphi_n\co  U\to V$ such that $\tilde\varphi_n$
descend to embeddings $\varphi_n\co U/\Gamma_{\mathrm{geo}}\to V/\Gamma_n$.
$$
\CD
U @>\tilde\varphi_n>> V \\
@VpVV @VVp_nV \\
U/\Gamma_{\mathrm{geo}} @>>\varphi_n> V/\Gamma_n\\
\endCD
$$
The embeddings $\tilde\varphi_n$ converge $C^r$--uniformly to
the inclusion.
\end{thm}

\begin{proof}[Sketch of the Proof]
Let $K=\{g\in G\co  g(\overline{V})\cap\overline{V}\ne\emptyset\}$, so
$K$ is compact. By~\ref{lok: discrete groups are open}, 
if we take $n$ sufficiently large, then
for any $\gamma\in\Gamma_{\mathrm{geo}}\cap K$ there exists a unique
$\gamma_n\in \Gamma_n$ that is close to $\gamma$.
It defines an injective map $r_n$ of $\Gamma_{\mathrm{geo}}\cap K$
into $\Gamma_n$.

Our goal is to construct embeddings $\varphi_n\co U\to V$ that are
close to the inclusion
and $r_n$--equivariant 
(in the sense that $\varphi_n(\gamma(x))=r_n(\gamma)\varphi_n(x)$ 
for $\gamma\in \Gamma_{\mathrm{geo}}$).
The construction is non-trivial and we refer to~\cite{Lok} for a complete proof.
\end{proof}

\begin{lem}\label{id component is compact} 
Let $G$ be a Lie group and
$\Gamma_n $ be a sequence of discrete groups that
converges geometrically to $\Gamma=\Gamma_{\mathrm{geo}}$.
If the identity component $\Gamma_0$ of $\Gamma_{\mathrm{geo}}$
is compact, the sequence $\{\Gamma_n \}$ has
unbounded torsion (ie for any positive integer $N$,
there is $\gamma\in \Gamma_{n(N)}$ of finite order
which is greater than $N$).
\end{lem}

\begin{proof} Choose $\epsilon>0$ so small that
the closed $4\epsilon$--neighborhood
$U$ of $\Gamma_0$ is disjoint from $\Gamma\setminus\Gamma_0$
(this is possible since $\Gamma_0$ is compact).
Given $M>1$, for large $n$, $\Gamma_n$
and $\Gamma$ are $\epsilon/M$--Hausdorff close on
the compact set $U$.
Take an arbitrary element $g\in \Gamma_0$
in $\epsilon$--neighborhood of the identity
$e$ and choose $g_n\in \Gamma_n$ so that it is 
$\epsilon/M$--close to $g$.

First, show that $g_n$ has finite order. Suppose not.
Since $U$ is compact and $\langle g_n\rangle$ is discrete and infinite,
there is a smallest $k>1$ with $(g_n)^k\notin U$.
So, $(g_n)^{k-1}\in U$.
The metric is left invariant, hence,
$$d((g_n)^{k-1},(g_n)^k)=
d(g_n,e)<d(g,e)+d(g_n,g)<\epsilon+\epsilon/M<2\epsilon.$$
Since $\Gamma_n$ and $\Gamma$ are $e/M$--Hausdorff close on
$U$ and since $\Gamma\cap U=\Gamma_0$,
we get $d((g_n)^{k-1},\Gamma_0)<\epsilon/M<\epsilon$.
This implies $d((g_n)^k,\Gamma_0)<3\epsilon$ and, thus, $(g_n)^k\in U$,
a contradiction. Thus $g_n$ has finite order.

We have just showed that any neighborhood of $g$ contains an element of
finite order $g_n\in \Gamma_n$ (for $n$ large enough).
So we get a sequence ${g_n}$ converging to $g$.
Assume the torsion of ${\Gamma_n}$ is uniformly bounded by $N$.
Then, for any finite order element $h\in \Gamma_n$, we have
$h^{N!}=e$.
In particular $e=(g_n)^{N!}$ converges to $g^{N!}$.
We, thus, proved that any element $g$ in $\Gamma_0$
that is $\epsilon$--close to
the identity satisfies $g^{N!}=e$.
This is absurd.
(Indeed, it would mean that the map
$\Gamma_0\to \Gamma_0$ that takes $g$ to $g^{N!}$
maps the $\epsilon$--neighborhood of the identity to the identity.
But the map is clearly a diffeomorphism on some
neighborhood of the identity).
Thus, ${\Gamma_n}$ has unbounded torsion as desired.
\end{proof}

\begin{lem}\label{id component is nilpotent} 
Let $G$ be a Lie group and let
$\Gamma_n $ be a sequence of discrete groups that
converges geometrically to $\Gamma=\Gamma_{\mathrm{geo}}$.
Then the identity component $\Gamma_0$ of $\Gamma_{\mathrm{geo}}$
is nilpotent.
\end{lem}

\begin{proof} Denote the identity
component by $G_0$. The group $\Gamma_{\mathrm{geo}}$
is closed, so it is a Lie group. 
So the identity component $\Gamma_0$ of $\Gamma_{\mathrm{geo}}$
is a Lie group;
we are to show that $\Gamma_0\le G_0$ is nilpotent.
Let U be a neighborhood of the identity that lies in a
Zassenhaus neighborhood in $G_0$ (see~\cite[8.16]{Rag}).
Being a Lie group $\Gamma_0$ is generated by any
neighborhood of the identity, in particular by $V=U\cap\Gamma_0$.
To show $\Gamma_0$ is nilpotent, it suffices to check that,
for some $k$,
$V^{(k)}=\{e\}$, that is, any iterated commutator of weight $k$
with entries in $V$ is trivial~\cite[8.17]{Rag}.
Fix an iterated commutator
$[v_1\dots[v_{k-2}[v_{k-1},v_k]]\dots]$ with $v_k\in V$
and choose a sequences $g_k^n\to v_k$ where $g_k^n\in \Gamma_n$.
For $n$ large enough, the elements $g_k^n,\dots g_k^n$
lie in a Zassenhaus neighborhood.
Then by Zassenhaus--Kazhdan--Margulis
lemma~\cite[8.16]{Rag} 
the group $\langle g_k^n,\dots g_k^n\rangle$ lies in a connected
nilpotent Lie subgroup of $G_0$.
The class of any connected nilpotent Lie subgroup of
$G_0$ is bounded by $\mathrm{dim}(G_0)$.
Thus, for $k>\mathrm{dim}(G_0)$, the $k$--iterated commutator
$[g^n_1\dots[g^n_{k-2}[g^n_{k-1},g^n_k]]\dots]$ is trivial,
for large $n$.
This implies
$[v_1\dots[v_{k-2}[v_{k-1},v_k]]\dots]$ is trivial and, 
therefore, $\Gamma_0$ is nilpotent.
\end{proof}

\section{Algebraic convergence}
\label{S:alg conv}
The set of representations $\mathrm{Hom}(\pi, G)$
of a group $\pi$ into a topological group $G$
can be given the so-called {\it algebraic}
topology (also called pointwise convergence topology).
Namely, a sequence of representations $\rho_k$ converges
to $\rho$ provided $\rho_k(\gamma)$ converges to $\rho(\pi)$ 
for each $\gamma\in\pi$. 
In this section we discuss how
algebraic and geometric convergences interact.

It follows from definitions that geometric limit
always contains the algebraic one. More precisely,
if $\rho_k$ converges
to $\rho$ algebraically and $\rho_k(\pi)$ converges to
$\Gamma_{\mathrm{geo}}$ geometrically, then
$\rho_k(\pi)\subset\Gamma_{\mathrm{geo}}$.
In particular, if $\Gamma_{\mathrm{geo}}$ is discrete,
so is $\rho(\pi)$. 

\begin{thm}\label{geometric limit discrete} 
Let $\pi$ be a finitely generated group
and let $X$ be a negatively curved
symmetric space.
Let $\rho_k\in\mathrm{Hom}(\pi,G)$ be a sequence of
representations converging algebraically to a representation $\rho$
where, for every $k$, the group $\rho_k(\pi)$ is discrete
and the sequence $\{\rho_k(\pi)\}$ has uniformly bounded torsion.
Suppose that $\rho(\pi)$ is an infinite group
without nontrivial normal nilpotent subgroups.

Then the closure of $\{\rho_k(\pi)\}$ in the geometric topology
consists of discrete groups. In particular,
the group $\rho(\pi)$ is discrete.
\end{thm}
\begin{proof} Passing to subsequence,
we can assume that $\{\rho_k(\pi)\}$ converges 
geometrically to $\Gamma_{\mathrm{geo}}$. 
Being a closed subgroup $\Gamma_{\mathrm{geo}}$ is a Lie group.
Suppose, arguing by contradiction, 
that the identity component $\Gamma_0$ of
$\Gamma_{\mathrm{geo}}$ is nontrivial.
By lemma~\ref{id component is nilpotent}, 
$\Gamma_0$ is a nilpotent Lie group.

Since $\rho(\pi)$ does not have a nontrivial normal nilpotent subgroup,
$\rho(\pi)\cap \Gamma_0=\{e\}$, so $\rho(\pi)$ is discrete.
Using the Selberg lemma we find a torsion-free subgroup
$\Pi\le\rho(\pi)$ of finite index.
The group $\Pi$ is infinite since $\rho(\pi)$ is.
Notice that the group
$\Pi$ may have at most two fixed points at infinity (any isometry
that fixes at least three points is elliptic~\cite[p 84]{BGS}).
The next goal is to show that $\Pi$ is virtually nilpotent.

According to~\cite[3.3.1]{Bow1},
the nilpotence of $\Gamma_0$ implies one
the following mutually exclusive conclusions:

{\it Case 1}\stdspace $\Gamma_0$ has a fixed point in $X$
(hence $\Gamma_0$ is compact which is is impossible 
by~\ref{id component is compact}).

{\it Case 2}\stdspace $\Gamma_0$ fixes a point
$p\in\partial X$ and acts freely
on $X\cup\partial X\setminus\{p\}$.

Since $\Pi$ normalizes $\Gamma_0$ (in fact,
$\Gamma_0$ is a normal subgroup of $\Gamma_{\mathrm{geo}}$)
it has to fix $p$ too. If $p$ is the only fixed point
of $\Pi$, then every nontrivial 
element of $\Pi$ is parabolic~\cite[8.9P]{EO}.
Any parabolic preserves horospheres centered at $p$~\cite[p 84]{BGS},
therefore, according to~\cite{Bow2} $\Pi$ is virtually nilpotent.

If $\Pi$  fixes two points at infinity, then
it acts freely and properly discontinuously on the geodesic
joining the points. Hence $\Pi\cong\Bbb Z$, the fundamental group
of a circle.

{\it Case 3}\stdspace $\Gamma_0$ has no fixed points in $X$
and preserves setwise some bi-infinite geodesic.
In this case the fixed point set of $\Gamma_0$ is the
endpoints $p$ and $q$ of the geodesic. Indeed, $\Gamma_0$ fixes each of
the endpoints, since $\Gamma_0$ is connected. Assume $\Gamma_0$
fixes some other point of $\partial X$.
Then any element of $\Gamma_0$ is elliptic (see~\cite[p 85]{BGS})
and, as such, it fixes the bi-infinite geodesic pointwise.
Thus $\Gamma_0$ has a fixed point in  $X$
which contradicts the assumptions of Case 3.

Since $\Pi$ normalizes $\Gamma_0$, it leaves the set $\{p,q\}$
invariant. Moreover, $\Pi$ preserves $\{p,q\}$ pointwise
because $\Pi$ contains no elliptics (any isometry that flips
$p$ and $q$ has a fixed point on the geodesic joining
$p$ and $q$).
Therefore, $\Pi\cong\Bbb Z$ as before.
Hence, $\Pi$ is virtually nilpotent.

Thus, $\rho(\pi)$ has a nilpotent subgroup of finite index
and, therefore, has a normal nilpotent subgroup of finite index.
This contradicts the assumption that $\rho(\pi)$
is an infinite group without nontrivial normal nilpotent subgroups.
\end{proof}

\begin{rmk}\label{remark on no normal nilpotent subgroup}
Let $\Gamma$ be a {\it torsion-free discrete} subgroup
of\/ $\mathrm{Isom}(X)$ that has a nontrivial
normal nilpotent subgroup.
Then $\Gamma$ is, in fact, virtually nilpotent.
Indeed, repeating the arguments above, we see that $\Gamma$ fixes 
a point at infinity and, hence, is virtually nilpotent.
Conversely, any virtually nilpotent group clearly has
a normal nilpotent subgroup of finite index.
\end{rmk}

\begin{lem} \label{algebraic limit is faithful}
Let $G$ be a Lie group and let $\pi$
be a finitely generated group
without nontrivial normal nilpotent subgroups.
Let $\rho_k\in\mathrm{Hom}(\pi,G)$ be a sequence of discrete
faithful representations that converges algebraically
to a representation $\rho$.
Then $\rho$ is faithful.
\end{lem}

\begin{proof} Denote the identity component of $G$ by $G_0$. 
It suffices to prove that the group 
$K=\mathrm{Ker}(\rho)$ is nilpotent.
Let $V$ be a set of generators for $K$ (maybe infinite).
To show $K$ is nilpotent, it suffices to check that,
for some $m$,
$V^{(m)}=\{e\}$, that is, any iterated commutator of weight $m$
with entries in $V$ is trivial~\cite[8.17]{Rag}.
Fix an iterated commutator
$[v_1\dots[v_{m-2}[v_{m-1},v_m]]\dots]$ with $v_m\in V$.

For $k$ large enough,
the elements $\rho_k(v_1),\dots ,\rho_k(v_m)$
lie in a Zassenhaus neighborhood of $G_0$.
Then by Zassenhaus--Kazhdan--Margulis
lemma~\cite[8.16]{Rag} the group
$\langle\rho_k(v_1),\dots ,\rho_k(v_m)\rangle$ lies in a connected
nilpotent Lie subgroup of $G_0$.
The class of any nilpotent subgroup of $G_0$ is
bounded by $\mathrm{dim}(G_0)$.
Thus, for $m>\mathrm{dim}(G_0)$, the $m$--iterated commutator
$$[\rho_k(v_1)\dots[\rho_k(v_{m-2})[\rho_k(v_{m-1}),\rho_k(v_m)]]\dots]$$
is trivial, for large $k$.
Since $\rho_k$ is faithful
$[v_1\dots[v_{m-2}[v_{m-1},v_m]]\dots]$ is trivial and we are done.
\end{proof}

\begin{cor} \label{injective, discrete is closed}
Let $X$ be a negatively curved symmetric space
and let $\pi$ be a finitely generated, torsion-free, discrete 
subgroup of\/ $\mathrm{Isom}(X)$ that is not virtually nilpotent.
Then the set of faithful discrete representations
of $\pi$ into $\mathrm{Isom}(X)$ is a closed subset of\/ 
$\mathrm{Hom}(\pi,\mathrm{Isom}(X))$.
\end{cor}
\begin{proof} 
Let a sequence $\rho_k$ of faithful discrete 
representations converge to $\rho\in\mathrm{Hom}(\pi,\mathrm{Isom}(X))$.
According to~\ref{remark on no normal nilpotent subgroup} 
the group $\pi$ has no normal nilpotent subgroup.
Then~\ref{algebraic limit is faithful} implies $\rho$ is faithful.
By compactness $\rho_k(\pi)$ has a subsequence
that converges geometrically to $\Gamma_{\mathrm{geo}}$
which is a discrete group by~\ref{geometric limit discrete}.
Therefore, $\rho(\pi)\subset\Gamma_{\mathrm{geo}}$
is also discrete as wanted.  
\end{proof}

\begin{cor} \label{geom limit of discrete injective is discrete}
Let $X$ be a negatively curved symmetric space
and let $\pi$ be a finitely generated, torsion-free, discrete 
subgroup of\/ $\mathrm{Isom}(X)$ that is not virtually nilpotent.
Suppose $\rho_k$ is a sequence of injective, discrete representations
of $\pi$ into $\mathrm{Isom}(X)$ that converges algebraically.
Then the closure of $\{\rho_k(\pi)\}$ in the geometric topology
consists of discrete groups.
\end{cor}
\begin{proof}
Combine~\ref{geometric limit discrete},~\ref{remark on no normal nilpotent subgroup},
and~\ref{injective, discrete is closed}.
\end{proof}

\section{A compactness theorem}
\label{S:comp thm}
In this section we state a compactness
theorem for the space of conjugacy classes
of faithful discrete representations of a group 
$\pi$ into the isometry group of a negatively curved 
symmetric space.
The proof follows from the work of Bestvina and Feighn~\cite{BF}
based on ideas of Rips 
(see the review of~\cite{BF}
by Paulin in MR$96h:20056$). 
Earlier versions of the theorem have
been proved by Thurston, Morgan and Shalen~\cite{Mor}.

Let $X$ be a negatively curved symmetric space
and $\mathrm{Isom}(X)$ be the isometry group of $X$.
We equip $\mathrm{Hom}(\pi,\mathrm{Isom}(X))$
with the algebraic topology.

\begin{thm}\label{compactness thm}
Let $X$ be a negatively curved symmetric space 
and let $\pi$ be a discrete, finitely presented subgroup 
of the isometry group of $X$. 
Suppose that $\pi$ is not virtually nilpotent and does not split 
as an HNN--extension or a nontrivial amalgamated product over a
virtually nilpotent group. 

Then the space of conjugacy classes of faithful discrete
representations of $\pi$ into $\mathrm{Isom}(X)$
is compact. 
\end{thm}

Here is an example of a group that does not
split over a virtually nilpotent group.

\begin{prop}\label{aspherical do not split}
Let $M$ be a closed aspherical $n$--manifold such that
any nilpotent subgroup of $\pi_1(M)$
has cohomological dimension $\le n-2$.
Then $\pi_1(M)$ does not split as a nontrivial
amalgamated product or HNN--extension over a
virtually nilpotent group.
\end{prop}

\begin{proof} 
Assume, by contradiction, that $\pi=\pi_1(M)$
is of the form $\Gamma_1\ast_N \Gamma_2$ or $\Gamma_0\ast_N$
where $N$ is a virtually nilpotent group
and $\Gamma_k$ is a proper subgroup of $\pi$, for $k=0,1,2$.

First, suppose that both $\Gamma_1$ and $\Gamma_2$
have infinite index in $\pi$. Note that by the definition
of HNN--extension, $\Gamma_0$ has infinite index in $\pi$.
Then it follows from the Mayer--Vietoris 
sequence~\cite[VII.9.1, VIII.2.2.4(c)]{Bro} that 
$$\mathrm{cd}(\pi)\le\underset{k}{\max}
(\mathrm{cd}(\Gamma_k), \mathrm{cd}(N)+1).$$
The cohomological dimension of $\pi=\pi_1(M)$ is $n$
because $M$ is a closed aspherical $n$--manifold.
Since $|\pi \co \Gamma_k|=\infty$, $\Gamma_k$ is the fundamental
group of a noncompact manifold of dimension $n$, hence
$\mathrm{cd}(\Gamma_k)<\mathrm{cd}(\pi)$. 
Finally, $\mathrm{cd}(N)\le n-2$ and we get a contradiction.

Second, assume that $\pi=\Gamma_1\ast_N \Gamma_2$ 
and, say, $\Gamma_1$ has finite index in $\pi$.
Look at the map $i_*\co H_n(\Gamma_1)\to H_n(\pi)$
induced by the inclusion $i\co \Gamma_1\to\pi$
in the homology with $w_1(M)$--twisted integer coefficients.
It is proved in~\cite[III.9.5(b)]{Bro} 
that the image of $i_*$
is generated by $|\pi \co \Gamma_1|\cdot[M]$
where $[M]$ is the fundamental class of $M$.

Since $\Gamma_1$ has finite index in $\pi$,
the group $N$ is of finite index in $\Gamma_2$ because
$N=\Gamma_1\cap\Gamma_2$. 
In particular, 
$\mathrm{cd}(\Gamma_2)=\mathrm{cd}(N)\le n-2$~\cite[VIII.3.1]{Bro}.

Look at the $n$th term
of the Mayer--Vietoris sequence with 
$w_1(M)$--twisted integer coefficients~\cite[VII.9.1]{Bro}:
\def\a{\longrightarrow}
$$\cdots\a H_n N\a H_n\Gamma_1\oplus H_n\Gamma_2\a  H_n\Gamma
\a H_{n-1}N\a \cdots$$
Using the cohomological dimension assumption on $N$,
we get $0=H_n(N)=H_{n-1}(N)=H_n(\Gamma_2)$,
and hence the inclusion $i\co \Gamma_1\to\pi$
induces an isomorphism of $n$th homology.

Thus $|\pi \co \Gamma_1|=1$ which contradicts
the assumption the $\Gamma_1$ is a proper subgroup of $\pi$.  
\end{proof}

\begin{cor} \label{word-hyperbolic dimension at least three}
Let $M$ be a closed aspherical
manifold of dimension $\ge 3$
with word-hyperbolic fundamental group.
Then $\pi_1(M)$ does not split as a nontrivial
amalgamated product or HNN--extension over a
virtually nilpotent group.
\end{cor}

\begin{proof} 
It is well known that any 
virtually nilpotent subgroup of a
word-hyperbolic group is 
virtually cyclic.
Since any virtually cyclic group has cohomological dimension
one, proposition~\ref{aspherical do not split} applies. 
\end{proof}

\section{The Main theorem}
\label{S:main thm}
Throughout this section 
$\iota$ is an invariant of continuous maps of $K$ into
oriented manifolds. 
We also use the letter $\iota$ for the corresponding invariant of 
discrete representations. 
Let $\pi$ be the fundamental group of a
finite CW--complex $K$.

\begin{thm} \label{main theorem}
Let $X$ be a contractible homogeneous Riemannian manifold
and $G$ be a transitive group of orientation-preserving
isometries.
Let $\rho_k\in\mathrm{Hom}(\pi, G)$
be a sequence of representations such that 
\begin{itemize}
\item the groups $\rho_k(\pi)$
are torsion-free, and
\item $\rho_k$ converges 
algebraically to $\rho$, and
\item $\rho_k(\pi)$ converges geometrically to 
a discrete group $\Gamma_{\mathrm{geo}}$. 
\end{itemize}
Let $f\co K\to X/\Gamma_{\mathrm{geo}}$ be a continuous map that
induces the homomorphism 
$\rho\co  \pi\to\rho(\pi)\subset\Gamma_{\mathrm{geo}}$.
Then $\iota (\rho_k)=\iota (f)$ for all large $k$. 
\end{thm}

\begin{proof} According to~\ref{lok: discrete groups are open}, 
we can assume that the groups $\rho_k(\pi)$ are discrete.
Being a limit of torsion-free groups the
discrete group $\Gamma_{\mathrm{geo}}$
is itself torsion-free by~\ref{limit is torsion free}.
Thus $\Gamma_{\mathrm{geo}}$ acts freely
and properly discontinuously on $X$,
so $X/\Gamma_{\mathrm{geo}}$ is a manifold.

We consider the universal covering  $q\co \tilde{K}\to K$.
Since $K$ is a finite complex,
one can choose a finite connected subcomplex $D\subset \tilde{K}$
with $q(D)=K$.
(Pick a representing cell in every orbit,
the union of the cells is a finite subcomplex that is mapped onto
$Y$ by $q$.
Adding finitely many cells, one can assume the complex
is connected.)

According to section \ref{S:Invariants of representations},
the representation $\rho$ defines a continuous 
$\rho$--equivariant map 
$\tilde f\co \tilde{K}\to X$ 
which is unique up to $\rho$--equivariant homotopy.
The map $\tilde{f}$ descends to a continuous map
$\bar f\co K\to X/\rho(\pi)$.

The identity map $\mathrm{id}\co X\to X$ is equivariant with
respect to the inclusion $\rho(\pi)\hookrightarrow\Gamma_{\mathrm{geo}}$,
therefore, $\tilde{f}$ is equivariant with respect to
the homomorphism $\rho\co \pi\to\Gamma_{\mathrm{geo}}$.
We denote by $f$ the composition of $\bar f$ and the covering
$X/\rho(\pi)\to X/\Gamma_{\mathrm{geo}}=N$ induced
by the inclusion $\rho(\pi)\subset \Gamma_{\mathrm{geo}}$.

Let $U\subset X$ be an open relatively compact
neighborhood of $\tilde{f}(D)$.
We are in position to apply Theorem~\ref{lok: lipschitz convergence}.
Thus, if $k$ is large enough, there is a sequence of embeddings
$\tilde{\varphi}_k\co U\to V\subset X$
that converges to the inclusion and is
$r_k$--equivariant (recall that $r_k$ was defined in the 
proof of~\ref{lok: lipschitz convergence}).
So the map $\tilde{\varphi}_k\circ\tilde{f}\co D\to V$ is
$r_k\circ\rho$--equivariant.
By the very definition of $r_k$ we have $r_k\circ\rho=\rho_k$
whenever the left hand side is defined.
We now extend the map
$\tilde{\varphi}_k\circ\tilde{f}\co D\to V$
by equivariance to a $\rho_k$--equivariant
map $\tilde{\xi}_k\co \tilde{K}\to X$.
The map $\tilde{\xi}_k$ descends to
$\xi_k\co K\to V/\rho_k(\pi)\subset X/\rho_k(\pi)$.
Notice that by construction $\xi_k=\varphi_k\circ f$.
Since $\varphi_k$ converge to the inclusion, 
$\varphi_k$ is orientation-preserving for large $k$.
Thus $\iota(\xi_k)=\iota(f)$ for large $k$
because $\iota$ is an invariant of maps into oriented manifolds.

Notice that $\iota(\xi_k)=\iota(\rho_k)$ since the map $\xi_k$
is $\rho_k$--equivariant 
(according to section \ref{S:Invariants of representations}
any $\rho_k$--equivariant
map can be used to define $\iota(\rho_k)$).
Therefore, for large $k$,
$\iota(\rho_k)=\iota(f)$ as desired. 
\end{proof}

\begin{rmk} 
Let $X$ be a negatively curved symmetric space.
Then, according to~\ref{geometric limit discrete},
the assumption ``$\Gamma_{\mathrm{geo}}$ is discrete''
of the theorem~\ref{main theorem} can be replaced
by ``$\rho_k(\pi)$ are discrete and 
$\rho(\pi)$ is an infinite group without nontrivial
normal nilpotent subgroups''.
\end{rmk}

\begin{rmk} In some cases the conclusion of the 
theorem~\ref{main theorem} can be improved to 
``$\iota (\rho_k)=\iota(\rho)$ for all large $k$''. 
Since $\iota (\rho)=\iota (\bar f)$, 
it suffices to show that $\iota(\bar f)=\iota (f)$. 

This is clearly true when $\rho(\pi)=\Gamma_{\mathrm{geo}}$
in which case it is usually said 
that $\rho_k$ converges to $\rho$ strongly.
Another obvious example is when $\iota$ is a liftable invariant.
Finally, according to~\ref{liftable prop}, 
$\iota(\bar f)=\iota (f)$ provided
$f$ is homotopic to an embedding
(ie a homeomorphism onto its image).
\end{rmk}

\begin{cor}\label{cor of main thm}
Let $X$ be a negatively curved symmetric 
space and let $K$ be a finite cell complex
such that $\pi=\pi_1(K)$ is a torsion-free,
not virtually nilpotent,
discrete subgroup of\/ $\mathrm{Isom}(X)$. 
Let $\rho_k$ be a sequence of discrete injective
representations of $\pi$ into the group of orientation-preserving
isometries of $X$. 
Suppose that $\rho_k$ converges in the pointwise convergence topology.

Then $\iota(\rho_k)=\iota(\rho_{k+1})$ for all large $k$.
\end{cor}
\begin{proof} The result follows from~\ref{main theorem} 
and~\ref{geom limit of discrete injective is discrete}.
\end{proof}

\begin{rmk} The results of this section certainly hold if $\iota$ 
is any invariant of maps rather that
an invariant of maps into {\it oriented} manifolds.
For such a $\iota$ we do not have to assume that
the isometry groups preserve orientation.
\end{rmk}

\section{Hyperbolic manifolds up to intersection preserving
homotopy equivalence}
\label{S:intersection}
A homotopy equivalence $f\co N\to L$ of oriented $n$--manifolds is called
{\it intersection preserving} if, for any $m$ and any pair of (singular) 
homology classes
$\alpha\in H_m(N)$ and $\beta\in H_{n-m}(N)$, their intersection number in
$N$ is equal to the intersection number of $f_*\alpha$ and $f_*\beta$
in $L$. For example, any map that is homotopic to an orientation-preserving
homeomorphism is an intersection preserving homotopy equivalence~\cite[13.21]{Dol}.

\begin{prop}\label{nilpotent groups}
Let $N$ be an oriented hyperbolic manifold with
virtually nilpotent fundamental group. Then the
intersection number of any two homology
classes in $N$ is zero.
\end{prop}
\begin{proof}
It suffices to prove that $N$ is homeomorphic
to $\Bbb R\times Y$ for some space $Y$.
First note that any torsion-free, discrete, virtually nilpotent group $\Gamma$
acting on a hyperbolic space $X$ must have either one or two 
fixed points at infinity 
(see the proof of~\ref{geometric limit discrete}). 
If $\Gamma$ has only one fixed point,
$\Gamma$ is parabolic and, hence, it preserves
all horospheres at the fixed point.
Therefore, if $H$ is such a horosphere,
$X/\Gamma$ is homeomorphic to $\Bbb R\times H/\Gamma$.
If $\Gamma$ has two 
fixed points, $\Gamma$ preserves a bi-infinite geodesic.
Hence $X/\Gamma$ is the total space of a vector bundle over
a circle and the result easily follows.  
\end{proof}

\begin{thm} \label{alg precompact implies finiteness of intersections}
Let $\pi$ be the fundamental group of a finite 
aspherical cell complex and let $X$ be a negatively curved symmetric 
space. Let $\rho_k$ be a sequence of discrete injective
representations of $\pi$ into the group of orientation-preserving
isometries of $X$. 
Suppose that $\rho_k$ is precompact in the pointwise convergence topology.

Then the sequence of manifolds $X/\rho_k(\pi)$ falls into
finitely many intersection preserving homotopy equivalence classes.
\end{thm}
\begin{proof} 
According to~\ref{nilpotent groups}, we can assume that
$\pi_1(K)$ is not virtually nilpotent.
Argue by contradiction.
Pass to a subsequence so that no two of the manifolds $X/\rho_k(\pi)$
are intersection preserving homotopy equivalent
and so that $\rho_k$ converges algebraically.

Being the fundamental group of a finite 
aspherical cell complex, $\pi$ has finitely generated homology.
Choose a finite set of generators of $H_*(\pi)$.
Using~\ref{cor of main thm}
we pass to a subsequence so that
$I_{n,\alpha,\beta}(\rho_k)=I_{n,\alpha,\beta}(\rho_{k+1})$
for any pair of generators
$\alpha\in H_m(N)$ and $\beta\in H_{n-m}(N)$.
Hence, the homotopy equivalence that induces $\rho_{k+1}\circ (\rho_k)^{-1}$
is intersection preserving.
\end{proof}

\begin{thm} \label{bound on intersections}
Let $K$ be a finite connected cell complex
such that $\pi_1(K)$ does not split 
as an HNN--extension or a nontrivial amalgamated product over a
virtually nilpotent group. 
Then, given a nonnegative integer $n$ and homology classes 
$\alpha\in H_m(K)$
and $\beta\in H_{n-m}(K)$,
there exists $C>0$ such that,

{\rm(1)}\stdspace for any continuous map $f\co K\to N$ of $K$ into an
oriented hyperbolic $n$--manifold $N$ that 
induces an isomorphism of fundamental groups, and

{\rm(2)}\stdspace for any embedding $f\co K\to N$ of $K$ into an
oriented hyperbolic $n$--manifold 
$N$ that induces a monomorphism of fundamental groups,

the intersection number of $f_*\alpha$ and
$f_*\beta$ in $N$ satisfies
$|\langle f_*\alpha, f_*\beta\rangle |\le C$.
\end{thm}

\begin{proof} First, we prove $(1)$.
Arguing by contradiction, consider a sequence
of maps $f_k\co K\to N_k$ that induce $\pi_1$--isomorphisms
of $K$ and oriented hyperbolic $n$--manifolds $N_k$
and such that $I_{n,\alpha,\beta}(f_k)$ is an unbounded
sequence of integers.
Pass to a subsequence so that no two integers $I_{n,\alpha,\beta}(f_k)$
are the same. 
Each map $f_k$ induces an injective discrete
representation $\rho_k$ of $\pi_1(K)$ into the group of orientation-preserving
isometries of a negatively curved symmetric space;
in particular $I_{n,\alpha,\beta}(\rho_k)=I_{n,\alpha,\beta}(f_k)$. 
Since in every dimension there are no more than four
negatively curved symmetric spaces, we can pass to a
subsequence to ensure that all the symmetric spaces
where $\rho_k(\pi_1(K))$ acts are isometric.

According to~\ref{nilpotent groups}, we can assume that
$\pi_1(K)$ is not virtually nilpotent. 
By~\ref{compactness thm} there exists a sequence of
orientation-preserving isometries $\phi_k$ of $X$
such that the sequence $\phi_k\circ\rho_k\circ(\phi_k)^{-1}$
is algebraically precompact. 
As we explained in the section~\ref{S:Invariants of representations},  
$I_{n,\alpha,\beta}(\rho_k)=
I_{n,\alpha,\beta}(\phi_k\circ\rho_k\circ(\phi_k)^{-1})$.  
Thus~\ref{cor of main thm} provides a contradiction.
Finally, according to~\ref{liftable prop}, $(1)$ implies $(2)$.
\end{proof}

\begin{cor} 
Let $\pi$ be the fundamental group of a finite 
aspherical cell complex that does not split 
as an HNN--extension or a nontrivial amalgamated product over a
virtually nilpotent group. 

Then, for any $n$, the class of hyperbolic $n$--manifolds with fundamental
group isomorphic to $\pi$ breaks into
finitely many intersection preserving homotopy equivalence classes.
\end{cor}
\begin{proof}
By~\ref{nilpotent groups}, we can assume that
$\pi_1(K)$ is not virtually nilpotent.
Thus the result follows 
from~\ref{alg precompact implies finiteness of intersections}
and~\ref{compactness thm}.
\end{proof}

\section{Vector bundles with hyperbolic total spaces}
\label{S:vector}
I am most grateful
to Jonathan Rosenberg for explaining to me the following nice fact.

\begin{prop}\label{P:classifying vector bundles}
Let $K$ be a finite CW--complex and $m$ be a positive integer.
Then the set of isomorphism classes of oriented real 
(complex, respectively) rank $m$ vector bundles over $K$
with the same rational Pontrjagin classes and the rational Euler class 
(rational Chern classes, respectively) is finite. 
\end{prop}

\begin{proof} To simplify notation 
we only give a proof for complex vector bundles
and then indicate necessary modifications for real 
vector bundles. 
We need to show that the ``Chern classes map'' 
$(c_1,\dots, c_m)\co [K, BU(m)]\to H^*(K,\Bbb Q)$
is finite-to-one. First, notice that $c_1$ classifies
line bundles~\cite[I.3.8, I.4.3.1]{Hir}, so we can assume $m>1$.

The integral $i$th Chern class 
$c_i\in H^*(BU(m),\Bbb Z)\cong [BU(m), K(\Bbb Z, 4i)]$ can be
represented by a continuous map $f_i\co BU(m)\to K(\Bbb Z, 2i)$
such that $c_i=f_i^*(\alpha_{2i})$ where $\alpha_{2i}$
is the fundamental class of $K(\Bbb Z, 2i)$.
(Recall that, by the Hurewicz theorem $H_{n-1}(K(\Bbb Z, n),\Bbb Z)=0$
and $H_n(K(\Bbb Z, n),\Bbb Z)\cong\Bbb Z$; the class
in $H^n(K(\Bbb Z, n),\Bbb Z)\cong\mathrm{Hom}(H_n(K(\Bbb Z, n),\Bbb Z);\Bbb Z)$
corresponding to the identity homomorphism is called the fundamental class
and is denoted $\alpha_n$.)

It defines ``Chern classes map'' 
$$c=(f_1,\dots, f_m)\co BU(m)\to 
K(\Bbb Z, 2)\times\dots\times 
K(\Bbb Z, 2m).$$ 
We now check that
the map induces an isomorphism on
rational cohomology. 
According to~\cite[7.5, 7.6]{GM} for even $n$, 
$H^n(K(\Bbb Z, n),\Bbb Q)\cong \Bbb Q[\alpha_n]$.  
By the K\"unneth formula 
$$H^*(\times_{i=1}^m K(\Bbb Z, 2i),\Bbb Q)\cong
\otimes_{i=1}^m H^*(K(\Bbb Z, 2i),\Bbb Q)\cong
\otimes_{i=1}^m \Bbb Q[\alpha_{2i}]\cong
\Bbb Q[\alpha_2,\dots \alpha_{2m}]$$
and under this ring isomorphism 
$\alpha_2^{s_1}\times\dots\times\alpha_{2m}^{s_m}$
corresponds to $\alpha_2^{s_1}\dots\alpha_{2m}^{s_m}$.
It is well known that 
$H^*(BU(m),\Bbb Q)\cong \Bbb Q[c_1,\dots, c_m]$~\cite[14.5]{MS}.
Thus the homomorphism 
$c^*\co H^*(\times_{i=1}^m K(\Bbb Z, 2i),\Bbb Q))\to H^*(BU(m),\Bbb Q)$
defines a homomorphism $\Bbb Q[\alpha_2,\dots \alpha_{2m}] \cong 
\Bbb Q[c_1,\dots, c_m]$ that takes $\alpha_{2i}$ to $c_i$.
Since this is an isomorphism, so is $c^*$ as promised.

Therefore $c$ induces an isomorphism on rational homology.
Then, since $BU(m)$ is simply 
connected, the map $c$ must be a rational homotopy
equivalence~\cite[7.7]{GM}. 
In other words the homotopy theoretic fiber $F_c$
of the map $c$ has finite homotopy groups.

Consider an oriented rank $m$ vector bundle
$f\co K\to BU(m)$ with characteristic classes
$c\circ f\co K\to\times_{i=1}^m K(\Bbb Z, 2i)$.
Our goal is to show that the map $c\circ f$ 
has at most finitely many nonhomotopic liftings to $BU(m)$.
Look at the set of liftings of 
$c\circ f$ to $BU(m)$ and try to construct homotopies
skeleton by skeleton using the obstruction theory. 
The obstructions lie in the groups of cellular
cochains of $K$ with coefficients in the
homotopy groups of $F_c$. 
(Note that the fibration $F_c\to BU(m)\to \times_{i=1}^m K(\Bbb Z, 2i)$ 
has simply connected base and fiber (since $m>1$),
so the coefficients are not twisted.)
Since the homotopy
groups of $F_c$ are finite, there are at most finitely many
nonhomotopic liftings.
This completes the proof for complex
vector bundles.

For oriented real vector 
bundles of odd rank 
the same argument works with 
$c=(p_1,\dots, p_{[m/2]})$, where $p_i$ is the 
$i$th Pontrjagin class.
Similarly, for oriented real vector 
bundles of even rank we set
$c=(e, p_1,\dots, p_{m/2-1})$, where $e$ is the Euler class.
(The case $m=2$ can be treated separately: 
since $SO(2)\cong U(1)$ any oriented rank two vector bundle
has a structure of a complex line bundle
with $e$ corresponding to $c_1$.
Thus, according to~\cite[I.3.8, I.4.3.1]{Hir}, 
oriented rank two vector bundles are in one-to-one
correspondence with $H^2(K,\Bbb Z)$.)
\end{proof}

\begin{rmk} A similar result is probably true for
nonorientable bundles. However, the argument given above
fails due to the fact $BO(m)$ is not simply connected 
(ie the map $c=(p_1,\dots, p_{[m/2]})$
of $BO(m)$ to the product of Eilenberg--MacLane spaces
is {\it not} a rational homotopy equivalence
even though it induces an isomorphism on
rational cohomology). 
For simplicity,
we only deal with oriented bundles.
\end{rmk}

\begin{rmk} \label{pontrjagin classifies}
Note that, if either
$m$ is odd or $m>\mathrm{dim}(K)$, the rational Euler class is
zero, and, hence rational Pontrjagin classes determine
an oriented vector bundle up to a finite number of possibilities.
\end{rmk}

\begin{cor} \label{normal bundles of immersions}
Let $M$ be a closed smooth manifold
and let $m$ be an integer such that 
either $m$ is odd or $m>\mathrm{dim}(M)$. 
Let $f_k\co B\to N_k$ be a sequence of smooth immersions 
of $B$ into complete locally symmetric nonpositively curved 
Riemannian manifolds $N_k$ with orientable normal bundle
and $\dim(N_k)=m+\dim(M)$.

Then the set of 
the normal bundles $\nu (f_k)$ of the immersions
falls into finitely many isomorphism classes.
\end{cor}
\begin{proof}
We can assume that all the manifolds $N_k$ are quotients
of the same symmetric space $X$ since in every dimension
there exist only finitely many symmetric spaces.

According to section~\ref{S:Invariants of representations}, 
the immersions $f_k$ induce representations 
$\rho_k\co \pi_1(M)\to\mathrm{Isom}(X)$ such that $\tau(f_k)=\tau(\rho_k)$.
The sequence $\tau(\rho_k)$ of vector bundles
breaks into finitely many isomorphism classes because
the representation variety $\mathrm{Hom}(\pi_1(M),\mathrm{Isom}(X))$
has finitely many connected components 
(see~\ref{thm on finitely many tan hom types}). 
In particular, there are only finitely many possibilities for the
total Pontrjagin class of $\tau(\rho_k)$.

The normal bundle of the immersion $f_k$ satisfies
$\nu (f_k)\oplus TM\cong\tau(f_k)$. Applying the total
Pontrjagin class, we get $p(\nu (f_k))\cup p(TM)=p(\tau(f_k))$.
The total Pontrjagin class of any bundle is a unit,
hence we can solve for $p(\nu (f_k))$.
Thus, there are only finitely many possibilities for $p(\nu (f_k))$.
Finally,~\ref{P:classifying vector bundles} and~\ref{pontrjagin classifies} 
imply that there are only finitely many possibilities for
$\nu(f_k)$.
\end{proof}

\begin{thm}
Let $M$ be a closed negatively curved
manifold of dimension $\ge 3$ and let $n>\dim(M)$ be an integer. 
Let $f_k\co B\to N_k$ be a sequence of smooth 
embeddings of $M$ into hyperbolic $n$--manifolds
such that for each $k$ 
\begin{itemize}
\item $f_k$ induces a monomorphism of fundamental groups, and
\item the normal bundle $\nu (f_k)$ of the embedding $f_k$
is orientable.
\end{itemize}
Then the set of 
the normal bundles $\nu (f_k)$  
falls into finitely many isomorphism classes.
In particular, up to diffeomorphism, 
only finitely many hyperbolic $n$--manifolds are
total spaces of orientable vector bundles over $M$.
\end{thm}
\begin{proof} 
Passing to covers corresponding to $f_{k*}$,
we can assume that $f_k\co M\to N_k$ 
induce isomorphisms of fundamental groups.
Arguing by contradiction, assume that $\nu_k=\nu (f_k)$ are
pairwise nonisomorphic.

Arguing as in the proof of~\ref{normal bundles of immersions},
we deduce that there are only finitely many possibilities for
the total Pontrjagin classes of $\nu_k$. 
Thus, according to~\ref{P:classifying vector bundles},
we can pass to a subsequence so that the (rational) Euler classes 
of $\nu_k$ are all different.
Denote the integral Euler class by $e(\nu_k)$.

First, assume that $M$ is orientable. 
Recall that, by definition, the Euler class
$e(\nu_k)$ is the image of the 
Thom class $\tau(\nu_k)\in H^m(N_k,N_k\backslash f_k(M))$
under the map $f_{k}^*\co H^m(N_k,N_k\backslash f_k(M))\to H^m(M)$.
According
to~\cite[VIII.11.18]{Dol} the Thom class   
has the property 
$\tau(\nu_k)\cap [N_k,N_k\backslash f_k(M)]=f_{k*}[M]$
where $[N_k,N_k\backslash f_k(M)]$ is 
the fundamental class of the pair 
$(N_k,N_k\backslash f_k(M))$ and $[M]$ is
the fundamental class of $M$.

Therefore, for any $\alpha\in H_m(M)$,
the intersection number of
$f_{k*}\alpha$ and $f_{k*}[M]$ in $N_k$ satisfies
$$I(f_{k*}[M], f_{k*}\alpha)=\langle\tau(\nu_k), f_{k*}\alpha\rangle=
\langle f^*\tau(\nu_k), \alpha\rangle=\langle e(\nu_k),\alpha\rangle.$$

Since $M$ is compact, $H_m(M)$ is finitely generated;
we fix a finite set of generators.
The (rational) Euler classes are all different, hence
the homomorphisms $\langle e(\nu_k),-\rangle\in\mathrm{Hom}(H_m(M),\Bbb Z)$ 
are all different. 
Then there exists a generator
$\alpha\in H_m(M)$ such that 
$\{\langle e(\nu_k),\alpha\rangle\}$ is an infinite set of integers.
Hence $\{I(f_{k*}[M], f_{k*}\alpha)\}$ is an infinite set
of integers. Combining~\ref{bound on intersections} 
and~\ref{word-hyperbolic dimension at least three}, 
we get a contradiction.

Assume now that $M$ is nonorientable.
Let $q\co\tilde M\to M$ be the orientable
two-fold cover. 
Any finite cover of aspherical manifolds
induces an injection on rational
cohomology~\cite[III.9.5(b)]{Bro}.
Hence $e(q^\#\nu_k)=q^*e(\nu_k)$
implies that the rational Euler classes of 
the pullback bundles $q^\#\nu_k$ are all different, and
there are only finitely many possibilities for
the total Pontrjagin classes of $q^\#\nu_k$.
Furthermore, the bundle map $q^\#\nu_k\to\nu_k$
induces a smooth two-fold cover of the total spaces,
thus the total space of $q^\#\nu_k$ is hyperbolic.
Finally, $\tilde M$ is a closed orientable
negatively curved manifold.
Thus, we get a contradiction 
as in the oriented case.
\end{proof}
\end{document}